\documentclass[10pt]{article}
\usepackage[mathscr]{eucal}
\usepackage{bbm}
\usepackage{dsfont}
\usepackage{color}
\usepackage{xcolor}
\usepackage{mathrsfs}
\usepackage{amssymb}
\usepackage{calligra}
\usepackage{fontenc}
\usepackage{amsbsy}
\usepackage{amsmath}
\usepackage{calligra}
\usepackage{fontenc}
\usepackage{amsfonts}
\usepackage[paperheight=10.75in, paperwidth=8.25in, margin=1in]{geometry}
\begin{document}
\input gpgmathprop31.sty
\pagenumbering{arabic}
\renewcommand{\thefootnote}{(\arabic{footnote})}
\title{On Local Regularity of Distributional Solutions\\ to the Navier--Stokes Equations} 
\author{ Giovanni P. Galdi
\thanks{Department of Mechanical Engineering and Materials Science, University of Pittsburgh, PA 15261}}
\date{}
\maketitle
\begin{abstract} We provide a sharp result that guarantees that a distributional solution satisfying the Prodi-Serrin condition is regular in the spatial variables. The solution does not need to belong to the (local) Leray-Hopf class.
\end{abstract}
\renewcommand{\theequation}{1.\arabic{equation}}
\setcounter{section}{0}
\section{Introduction}
As is well known, Leray-Hopf solutions occupy a prominent place in the mathematical theory of the Navier-Stokes equations. This is because they exist on an arbitrary time interval $(0,T)$, without restrictions on the "size" of the data $\bfu_0$ \cite{Hopf,Leray}. More specifically, in the case of Cauchy problem, they are   
characterized by a  vector field $\bfu$ in the {\em Leray-Hopf} regularity class~\footnote{$L_\sigma^2(\real^3)$ is the subspace of $L^2(\real^3)$ of solenoidal functions. Other notations are standard, like $L^q$ and $W^{m,q}$ for Lebesgue and  Sobolev spaces, with corresponding norms $\|\cdot\|_{q}$ and $\|\cdot\|_{m,q}$, respectively, $L^r(I;X)$, $I$  real interval, $X$ Banach space, for Bochner spaces, etc. Other notations will be introduced at the beginning of Section \ref{PrepRes}.}
\be
\mathscr {L H}_T:=L^\infty(0,T;L^2_\sigma(\real^3))\cap L^2(0,T;W^{1,2}(\real^3))\,,
\eeq{1.1}
that solves the following equation
\be
\Int0T\Int{\real^3}{}\big(\bfu\cdot\partial_t\bfphi+ \bfu\cdot\Delta\bfphi+\bfu\cdot\nabla\bfphi\cdot\bfu\big)=-\Int{\real^3}{}\bfu_0\cdot\bfphi(0)
\eeq{1.2}
for some $\bfu_0\in L^2_\sigma(\real^3)$ and all $\bfphi\in \mathcal D_T:=\{\bfphi\in C_0^\infty([0,T)\times\real^3):\,\Div\bfphi=0\}$.
\par
However, it is equally known that, to date, fundamental problems such as their uniqueness, regularity and fulfillment of the energy balance (``energy equality") remain open.\footnote{We consider here the case of zero external forces. In this regard, it must be emphasized that if the initial data vanish and the force is chosen suitably, the uniqueness property has been disproved by the sharp counterexample shown in  \cite{AlCo}.}
\par
Since Leray's pioneering  contribution \cite[p. 227 and \S\,22]{Leray}, mathematicians have wondered what {\em further} regularity conditions a solution $\bfu\in \mathscr L\mathscr H_T$ must possess in order to satisfy the above properties. For example, if $\bfu\in\mathscr L\mathscr H_T$  and {\em in addition} satisfies
\be
\bfu\in L^{q'}(0,T; L^q(\real^3))\,,\ \ \frac2{q^\prime}+\frac3q=1\,,\ \ q>3\,,
\eeq{1.3}
then $\bfu$ is in $C^\infty((0,T)\times\real^3)$ \cite{GaMa,Giga,Sohr}, as well as  it is the only one \cite{KoSo,Ser1} among those solutions $\bfv\in \mathscr L\mathscr H_T$ corresponding to the same $\bfu_0$ and   satisfying the ``energy inequality"
\be  
\|\bfv(t)\|_2^2+2\int_0^t\|\nabla\bfv(s)\|_2^2{\rm d}s\le \|\bfu_0\|_2^2\,,\ \ t\in [0,T]\,.
\eeq{1.4}  
Furthermore, if $\bfu\in \mathscr L\mathscr H_T$ is such that
\be
\bfu\in L^{r}(0,T; L^s(\real^3))\,,\ \ \frac2{r}+\frac2s=1\,,\ \ s\ge 4\,,
\eeq{1.5}
then $\bfu$ satisfies the ``energy equality," that is, \Eqref{1.4} with ``$\le$" replaced by "$=$" and $\bfv\equiv\bfu$  \cite{Lions,Shin}.   
\par
Recently, the present author has started to investigate whether the condition that $\bfu$ belongs to the Leray-Hopf class, namely that $\bfu$ possesses {\em finite energy}, is indeed a fundamental prerequisite for the validity of all the regularity results above. The answer is that such an assumption is, in general, {\em redundant}, which makes, in a sense,  Leray-Hopf class not so special. Precisely,  suppose that $\bfu$ is {\em only} in $L^2_{\rm loc}([0,T)\times\real^3)$,  satisfies \Eqref{1.2} for all $\bfphi\in\cald_T$ and some $\bfu_0\in L^2_\sigma(\real^3)$, and is divergence-free. 
In \cite{GaP1} it is proven that if $\bfu$ satisfies \Eqref{1.5} with $r=s=4$, then $\bfu\in \mathscr L\mathscr H_T$ and, therefore,  obeys the energy equality. 
Furthermore, in \cite{GaP2} it is shown that if $\bfu$ obeys  \Eqref{1.3}, then automatically $\bfu\in \mathscr L\mathscr H_T$ and, therefore, $\bfu$ is, in particular, smooth and unique. Results in \cite{GaP1,GaP2} have been further extended and refined by several authors in different directions; see, e.g., \cite{BeCh,DiTa,Lai,Wu,Wu1}. 
\par
At this point, it seems natural to ask whether it is possible to prove a {\em local} version of the regularity results in \cite{GaP2}. In particular, whether the well-known Serrin regularity criterion \cite{Ser} -- refined by Struwe \cite{Stru} and Takahashi \cite{Tak} -- can also hold
without assuming that the solution $\bfu$ belongs to the local Leray-Hopf class. To state this question precisely, let $\mathscr O=B\times (t_1,t_2)$, with $B$ open ball, be an arbitrary space-time cylinder in $ \real^3\times (0,T)$ and suppose that $\bfu\in L^2(\mathscr O)$ is solenoidal and satisfies
\be
\Int{\mathscr O}{}\big(\bfu\cdot\partial_t\bfphi+ \bfu\cdot\Delta\bfphi+\bfu\cdot\nabla\bfphi\cdot\bfu\big)=0\,,\ \ \mbox{for all $\bfphi\in\cald(\mathscr O)$}\,,
\eeq{1.6}
where $\cald(\mathscr O):=\{\bfphi\in C^\infty_0(\mathscr O): \,\Div\bfphi=0\}.$ Let $L^{r',r}(\mathscr O)$ be the space of functions  with
$$
\|\bfu\|_{L^{r',r}(\mathscr O)}:=\left(\int_{t_1}^{t_2}\|\bfu(\cdot,t)\|_{r,B}^{r'}dt\right)^{1/r'}<\infty\,.
$$
The criterion then states that if $\bfu$ satisfies
\be
\bfu\in L^{q',q}(\mathscr O)\,,\ \ \frac2{q'}+\frac3{q}=1\,,\ \ q>3\,,
\eeq{1.8}
and {\em in addition} is locally in the Leray-Hopf class, namely, 
\be
\bfu\in L^{\infty,2}(\mathscr O),\ \nabla\bfu\in L^{2,2}(\mathscr O)\,,
\eeq{1.7}
then $\bfu$ is $C^\infty$ in the space variables and each derivative is $L^\infty$ in time on compact subdomains of $\mathscr O$.~\footnote{Actually, in \cite{Ser} the assumption \Eqref{1.7} is replaced by the seemingly weaker one $\bfu\in L^{\infty,2}(\mathscr O),\ \curl\bfu\in L^{2,2}(\mathscr O)$. However, it is not difficult to show that the latter and \Eqref{1.7} are, in fact, equivalent.}
The question we want to answer is therefore: {\em Is hypothesis \Eqref{1.7} really necessary?}
\par
Before addressing this question, we recall that,  as remarked in \cite{Ser}, equation  \Eqref{1.6} admits the potential-like solution $\hat{\bfu}=a\,\nabla\psi$ where $\psi=\psi(x)$ is harmonic in $B$ and $a=a(t)$ is {\em just} in $L^2(t_1,t_2)$.     
Therefore, for Serrin's result to be valid, it is essential that the hypotheses either include the $L^\infty$ {\em in time} property somewhere, or exclude potential-like solutions. 
\par
A first contribution to the question was furnished in   \cite{Sere}, for solutions satisfying \Eqref{1.8} and possessing an associated pressure field $p$, suitably defined. In place of \Eqref{1.7}, the author assumes that $\bfu,p$ and $\nabla\bfu$ are summable in certain $L^{r',r}$-spaces, with $r'<\infty$,  and shows $\bfu\in C(\mathscr O)$. The $L^\infty$-in time property mentioned above is hidden in the assumption on the pressure. Actually, the pressure field associated to $\hat{\bfu}$ is $\hat{p}(x,t)=-a'(t)\psi(x)-\half a^2(t)|\nabla\psi(x)|^2$, so that, recalling that $a$ is, at least, square integrable, requiring $\hat{p}\in L^{r',r}$ for {\em some} $r'\ge1$, implies $\hat{\bfu}\in L^{\infty,\infty}$.\par 
Another remarkable contribution has been recently given in \cite{DiTa}.\footnote{I am indebted to this paper for some of the methods used here.}  There, the authors replace \Eqref{1.7} with the assumption that $\bfu\in L^2([0,T];L_{uloc}^2(\real^3))$ is a mild solution to the Cauchy problem with initial data $\bfu_0\in L^2_\sigma(\real^3)$.  
This result, though interesting, is not exactly in the spirit of local regularity, in that it requires that $\bfu$ solves the initial-value problem. However, a requirement of this type is somehow necessary. In fact, since  boundedness in time for $\bfu$ is {\em not} assumed, one must rule out  solutions of the type $\hat{\bfu}$, and this is accomplished by the condition $\bfu_0\in L^2(\real^3)$.   
\par
Objective of this note is to furnish a rather complete answer to the question raised above. Precisely, let $\bfu\in L^2(\mathscr O)$ satisfy \Eqref{1.6}, and let $\bfu=\bfu_\sigma+\nabla\pi$ be its Helmholtz-Weyl decomposition in $B$. Then, in \theoref{3.1} we prove that if $\bfu_\sigma$ obeys \Eqref{1.8} and $\pi\in L^{\infty,1}(\mathscr O)$, it follows that $\bfu$ satisfies \Eqref{1.7} on every compact subdomain of $\mathscr O$ and, therefore, after Serrin and Struwe's results, $\bfu$ is smooth in the space variables.  Some remarks are in order. First, the assumption of $L^\infty$ in time is made only on the gradient part and cannot be removed, as the potential-like solution example shows. Secondly, since $\pi$ is harmonic and therefore $C^\infty$ in space, lack of spatial regularity can only come from $\bfu_\sigma$ and, in fact, assumption \Eqref{1.8} is needed only for $\bfu_\sigma$. Finally, the hypothesis on $\pi$ could be replaced by  $\bfu\in L^{\infty,s}(\mathscr O)$, for {\em some} $s>1$; see \remref{3.2}.   
\par
The proof of our theorem is rather straightforward and relies on a space-time regularity property for distributional solutions to homogeneous and non-homogeneous Stokes {\em systems} shown in Section 2; see \lemmref{HE} and \lemmref{Oseen}. Thanks to these results, in Section 3 we give a proof of our main finding (\theoref{3.1}), by first assuming $q\in[4,6]$ (\lemmref{3.1}), and then, also adopting  an elegant argument used in \cite{DiTa}, for all $q\in (3,4)\cup(6,\infty)$ (\lemmref{3.2}).   
We end this note with a remark about time regularity, suggesting that such a property could hold for $\bfu$, if it only holds for $\pi$.

\renewcommand{\theequation}{2.\arabic{equation}}

\section{Preparatory Results.}
\label{PrepRes}
We begin to recall some  notation. Open balls 
of $\real^3$ are indicated by the symbols $B$, $B'$ etc. We also set
$$
\langle\bfv,\bfw\rangle_{B}:=\int_B\bfu(x)\cdot\bfv(x){\rm d}x.
$$
A space-time cylinder in $\real^4_+:=\real^3\times (0,\infty)$ is the Cartesian product of an open ball in $\real^3$ with a bounded open interval $I\subset (0,\infty)$, and will be typically denoted  by $Q$, $Q'$,  etc. By $\cald(Q)$ we denote the space of solenoidal vector functions that are infinitely differentiable in $Q$ with compact support in $Q$.
 If $D$ is a  domain of $\real^3$ or $\real^4_+$, we will write $D'\Subset D$ whenever $D^\prime$ is a domain with $\bar{D'}\subset D$. For $q\in(1,\infty)$,  the following Helmholtz-Weyl decomposition of $L^q(B)$ holds
(see \cite[Chapter III]{Gab})
\be
L^q(B)=L^q_\sigma(B)\oplus G_q(B)\,.
\eeq{HelD}
where
$$\ba{ll}\medskip
L^q_\sigma(B):=\{\bfv\in L^q(B): \Div\bfv=0,\ \bfv\cdot\bfn|_{\partial B}=0\}\\
G_q(B):=\{\bfw\in L^q(B): \bfw=\nabla\pi\,,\ \pi\in W^{1,q}(B)\}\,.
\ea
$$ 
The projection operator from $L^q(B)$ onto $L^q_\sigma(B)$ is denoted by $\mathsf P$. It is known that $\mathsf P$ depends on $B$, but it is independent of $q$. For $\bfv\in L^q(B)$ we set 
$$
\bfv_\sigma:=\mathsf P\,\bfv\,,
$$
and $\mathsf Q:=\mathsf I-\mathsf P$, with $\mathsf I$ identity in $L^q$.
\par
\Bl Let $\bfv\in L^q(B)$, $1<q<\infty$, with $\Div\bfv=0$ in the sense of distributions, and let $\nabla\pi=\mathsf Q\,\bfv$. Then $\pi\in C^\infty(B)$. Moreover,  for any $B'\Subset B$, there is $C=C(B',B,k)>0$ such that 
\be
\sup_{x\in B'}|\partial_x^{|k|}\pi(x)|\le C\,\|\pi\|_{1,B}\,,\ \ \mbox{all $|k|\ge 0$}\,.
\eeq{ciacio}
\EL{har}
{\em Proof.} Since $\Div\bfv=0$, it follows that $\pi$ is harmonic in $B$ and, therefore, of class $C^\infty(B)$. Furthermore, from \cite[Remark V.3.2]{Gab} (see also \cite[\S 4]{Fuji}), we deduce that $\pi$ has the following representation
$$
\pi(x)=\int_BH^\delta(x-y)\pi(y){\rm d}y\,,\ \ x\in B'\,,
$$ 
where $H^\delta$ is infinitely differentiable and vanishes unless $|x-y|\in (\delta/2,\delta)$ with $\delta={\rm dist}\,\{\partial B,\partial B'\}$. Therefore, \Eqref{ciacio} follows.
\par\hfill$\square$\par
We shall next prove a similar regularity result for distributional solutions to the Stokes system in a space-time cylinder. This property could be deduced from \cite{OH}. However, in view of its relevance, we shall give here a full and independent proof. To this end,
 we recall the Oseen fundamental tensor solution, $
{I}\!\!\varGamma = \{\varGamma_{ij}\}$, to the Stokes equation \cite[\S\, VIII.3]{Gab}, defined as 
$$\ba{ll}\medskip
\varGamma_{ij}(x-y,t-\tau)=-\delta_{ij}\Delta\Psi(|x-y|,t-\tau)+\partial_{y_i}\partial_{y_j}\Psi(|x-y|,t-\tau)\,,\\
\Psi(r,s):=\left\{\ba{ll}\medskip\Frac{1}{4\pi^{3/2}s^{1/2}}\Frac{1}{r}\Int 0{r}{{\rm e}^{-\frac{\rho^2}{4s}}}{\rm d}\rho\ \ &\mbox{if $s>0$}\\
0\ \ &\mbox{if $s\le 0$}\ea\right.
\ea
$$
and satisfying
\be\left.\ba{ll}\medskip\partial_\tau\varGamma_{ij}+\Delta_y\varGamma_{ij}=0\\
\pde{\varGamma_{ij}}{y_j}=0\ea\right\}\ \ \mbox{$ \tau<t$}\,.
\eeq{GA}
\par
The following result holds.
\Bl Let $Q:=B\times (0,T)$, and let  $\bfv\in L^{1,q}(Q)$, $q>1$, with $\Div\bfv(t)=0$ in $B$, a.a. $t\in (0,T)$, in the sense of distributions. Then, if  
\be
\int_0^T\langle{\bfv},\partial_t\bfphi+\Delta\bfphi\rangle_B=0\ \ \mbox{for all $\bfphi\in\cald(Q)$}\,,
\eeq{C1}
it follows that $\bfv_\sigma\in C^\infty(Q')$ for any $Q':=B'\times(t_1^\prime,t_2^\prime)\Subset Q$. Moreover,  there is $C=C(Q',Q,\alpha,\beta)>0$ such that 
$$
\sup_{(x,t)\in Q'}|\partial^{\alpha}_t\partial^{|\beta|}_x \bfv_\sigma(x,t)|\le C\,\|\bfv_\sigma\|_{L^{1,1}(Q)}\ \ \ \mbox{for all $\alpha,|\beta|\ge 0$\,.}  
$$
\EL{HE}
{\em Proof.} In view of the decomposition \Eqref{HelD} and the assumption on $\bfv$, it follows that \Eqref{C1} is equivalent to
\be
\int_0^T\langle{\bfv_\sigma},\partial_t\bfphi+\Delta\bfphi\rangle_B=0\ \ \mbox{for all $\bfphi\in\cald(Q)$}\,.
\eeq{C2}
Let $Q^{\prime\prime}:=B^{\prime\prime}\times (t_1,t_2)$ be such that 
$${Q'}\Subset Q^{\prime\prime}\Subset  Q,$$ and  let  $\bfphi$ be an arbitrary element of $\cald(Q{''})$. For $\eta$ sufficiently small, the space-time (Friedrichs) mollifier, $\bfphi_\eta$, of $\bfphi$ belongs to $\cald(Q)$ and can  therefore be replaced in \Eqref{C2}. By a standard argument that uses the commutativity of mollification and differentiation operations, we then show that the mollifier, $\bfv_{\sigma\eta}$, of $\bfv_\sigma$ satisfies
\be\left.\ba{ll}\medskip
\partial_t\bfv_{\sigma \eta}-\Delta \bfv_{\sigma \eta}=\nabla p^{\eta}\\
\Div\bfv_{\sigma\eta}=0\ea\right\}
\ \ \mbox{in $B^{\prime\prime}\times (t_1,t_2)$}\,, 
\eeq{C3}
for some smooth, harmonic scalar function $p^{\eta}$.
Let $\delta=\dist\{\partial Q',\partial Q^{\prime\prime}\}$.
For $(x,t)\in Q'$, we denote by   $\varphi^\delta=\varphi^\delta(x-y,t-\tau)$ a smooth function that is equal to 1 if $|x-y|\le\delta/2$ and $|t-\tau|\le\delta/2$, while it is equal to 0 if either $|x-y|\ge\delta$ or $|t-\tau|\ge\delta$. Setting
$$
\varGamma_{ij}^\delta(x-y,t-\tau):=\varphi^\delta(x-y,t-\tau)\varGamma_{ij}(x-y,t-\tau)\,,
$$
from \Eqref{GA} and the properties of $\varphi^\delta$ we readily deduce that
\be\left.\ba{ll}\medskip\partial_\tau\varGamma_{ij}^\delta+\Delta_y\varGamma_{ij}^\delta=H_{ij}^\delta(x-y,t-\tau)\\
\pde{\varGamma_{ij}^\delta}{y_j}=0\ea\right\}\ \ \mbox{$ \tau<t$}\,.
\eeq{GAd}
where 
$$
H_{ij}^\delta=H_{ij}^\delta(z,s)
\ \ \mbox{is a $C^\infty$ function with $\supp (H_{ij}^\delta)\subset \{\delta/2<|z|+s<\delta\}$}\,.
$$ 
Set ${\bfvarGamma}_i^\delta=(\varGamma_{i1}^\delta,\varGamma_{i2}^\delta,\varGamma_{i3}^\delta)$,
and let
$$
\mathbb T(\bfu,p)=-p\,\mathbb I+\nabla\bfu+(\nabla\bfu)^\top
$$
with $\mathbb I$ identity matrix, 
denote the Cauchy stress tensor.
Then, for small $\varepsilon>0$ and $t\in(t_1+\varepsilon,t_2)$, we use the well-known Green's formula  
$$\ba{rl}\medskip
\Int{t_1}{t-\varepsilon}&\!\!\!\!\Int{B^{\prime\prime}}{}\{(\partial_\tau{\bfu_1}+\Delta \bfu_1-\nabla p_1)\cdot\bfu_2+
(\partial_\tau{\bfu_2}-\Delta \bfu_2+\nabla p_2)\cdot\bfu_1\}{\rm d}y\,{\rm d}\tau
\\
&=\Int{t_1}{t-\varepsilon}\Int{\partial B^{\prime\prime}}{}\big[\bfu_1\cdot\mathbb T(\bfu_2,p_2)-\bfu_2\cdot\mathbb T(\bfu_1,p_1)\big]\cdot\bfn {\rm d}\sigma_y+
\left.\Int{B^{\prime\prime}}{}\bfu_1(y,\tau)\cdot\bfu_2(y,\tau){\rm d}y\right|_{\tau=t_1}^{\tau=t-\varepsilon}\,,
\ea
$$
with $\bfu_1:={\bfvarGamma}_i^\delta(x-y,t-\tau)$, $i=1,2,3$, $p_1\equiv0$, and $\bfu_2:=\bfv_{\sigma\eta}(y,\tau)$, $p_2=p^\eta$. 
Employing \Eqref{C3} and \Eqref{GAd}, and observing that, by the properties of $\varphi^\delta$,
$$\ba{ll}\medskip
\partial_x^k\bfvarGamma^\delta_i(x-y,t-\tau)=\0\,,\ \ (x,t)\in Q'\,, \ \,(y,\tau)\in \partial B^{\prime\prime}\times (t_1,t-\varepsilon)\,,\ \     |k|\ge 0\,,\ i=1,2,3\,,\\ \bfvarGamma_i^\delta(x-y,t_1)=\0\,,\ \ (x,y)\in B'\times B^{\prime\prime}\,, \ \, i=1,2,3\,,
\ea
$$
we then show
\be
\int_{B^{\prime\prime}}(\bfv_{\sigma\eta})_j(y,t-\varepsilon)\varGamma_{ij}(x-y,\varepsilon)\varphi^\delta(x-y,\varepsilon){\rm d}y\,{\rm d}\tau=\int_{t_1}^{t-\varepsilon}\int_{B^{\prime\prime}}H_{ij}^\delta(x-y,t-\tau)(\bfv_{\sigma\eta})_j(y,\tau){\rm d}y\,{\rm d}\tau\,.
\eeq{1.2_n}
By \cite[\S\,5{\small 3}]{Oseen} (see also \cite[Lemma VIII.3.1]{Gab}) and the property of $\varphi^\delta$ we infer
\be\ba{rl}\medskip
\Lim{\varepsilon\to0}\Int{ B^{\prime\prime}}{}(\bfv_{\sigma\eta})_j(y,t-\varepsilon)&\!\!\!\!\varGamma_{ij}(x-y,\varepsilon)\varphi^\delta(x-y,\varepsilon){\rm d}y\,{\rm d}\tau\\
&\!\!\!\!=(\bfv_{\sigma\eta})_j(x,t)-\Frac1{4\pi}\Int{\partial B^{\prime\prime}}{}\Frac{x_j-y_j}{|x-y|}\varphi^\delta(y,t)\,\bfv_{\sigma\eta}(y,t)\cdot\bfn\,{\rm d}o_y=(\bfv_{\sigma\eta})_j(x,t)\,.\ea
\eeq{eps}
Thus, letting $\varepsilon\to0$ in \Eqref{1.2_n} and using \Eqref{eps}, and then letting $\eta\to 0$ in the resulting relation, we easily deduce (after redefining $\bfv_\sigma$ on a set of zero Lebesgue measure)
$$
v_{\sigma i}(x,t)=\int_{t_1}^t\int_{B^{\prime\prime}}H^\delta_{ij}(x-y,t-\tau)v_{\sigma j}(y,\tau){\rm d}y\,{\rm d}\tau\,,\ \ \ (x,t)\in Q'\,,
$$
from which the lemma follows.\par\hfill$\square$\par
\Bl
Let $\mathscr R_{t_1,t_2}:=\real^3\times (t_1,t_2)$, and let $\mathbb{F}=\{F_{ij}\}$ be a  sufficiently smooth tensor function in  (e.g., $\mathbb{F}$ is continuous in $\mathscr R_{t_1,t_2}$ and H\"older continuous in $x\in \real^3$ with its first derivatives, uniformly in $t$), with a bounded spatial support independent of $t$. Then, the vector function $\bfv$ with components
\be
v_i(x,t)=\int_{t_1}^t\int_{\real^3}\partial_{y_\ell}\varGamma_{ij}(x-y,t-\tau)\cdot F_{\ell j}(y,\tau)dy\,d\tau\,,\ \ (x,t)\in \mathscr R_{t_1,t_2}\,,\ \ i=1,2,3,
\eeq{SS}
is smooth as well (e.g., twice differentiable in the space variables and differentiable in time with all derivatives continuous in $\mathscr R_{t_1,t_2}$) and 
satisfies the equations
\be\left.\ba{ll}\medskip
\partial_t{\bfv}=\Delta\bfv+\nabla\Phi+\Div\mathbb F\\
\Div\bfv=0\ea\right\}\ \ \mbox{ in $\mathscr R_{t_1,t_2}$}
\eeq{SS1}
for a suitable smooth scalar field $\Phi$ (e.g., differentiable in space with each derivative continuous in $\mathscr R_{t_1,t_2}$). 
\par
Moreover, suppose
\be 
\mathbb F\in L^{\frac{q'}2,\frac q2}(\mathscr R_{t_1,t_2})\,,\ \ \frac2{q'}+\frac3q=1\,.
\eeq{BM1}
Then
$$
\bfv\in L^{r',r}(\mathscr R_{t_1,t_2})
$$
where
\be
\frac2{r'}+\frac3r=1\,,\  \ \ 2r\ge q>\frac{6r}{3+r}\,,
\eeq{BM2}
and 
\be
\|\bfv\|_{L^{r',r}}\le c\,\|\mathbb F\|_{L^{\frac{q'}2,\frac q2} }
\eeq{vf}
\EL{Oseen}
{\em Proof.} The first part of the lemma is a classical result of Oseen (see \cite[\S\S 5-7]{Oseen}).
To show the second part, we begin to recall the following estimate of the fundamental tensor, also due to Oseen \cite[p. 70]{Oseen},
\be
|\partial_{\xi_j}{I}\!\!\varGamma(\xi,s)|\le \Frac{C}{(|\xi|^2+s)^2}\,,\ \ \ j=1,2,3,
\eeq{SS_0}
for some constant $C>0$. Thus, taking into account \Eqref{SS_0}, we apply Young's inequality for convolution in \Eqref{SS} to get 
\be
\|\bfv(t)\|_{r}\le \int_{t_1}^{t}(t-\tau)^{-2+\frac{3}{2p}}\|\mathbb F(\tau)\|_{\frac q2}\,{\rm d}\tau\,,
\eeq{FN_0} 
where
\be
\frac 1r=\frac 1p+\frac 2q-1\,;\ \ \ \frac1r\le\frac2q\,.
\eeq{FN}
Next, we observe that the function $h(s)=s^{-\alpha}$, $\alpha:=2-\frac3{2p}$, is in the weak space $L_w^{1/\alpha}(0,T)$, $T>0$ for all $\alpha>0$. In particular, if $\alpha<1$, namely, 
\be
p\in[1,\mbox{$\frac32$})\,.
\eeq{fn}
we may use in \Eqref{FN_0}  the generalized Young inequality \cite[\S IX.4]{ReSi}
and obtain
$$
\|\bfv\|_{L^{r',r}}\le c_1\|\mathbb F\|_{L^{\frac{q'}2,\frac q2} }
$$
where
\be
\frac1{r'}=\alpha+\frac2{q'}-1=1+\frac2{q'}-\frac3{2p}\,.
\eeq{FN2}
Therefore, from \Eqref{BM1}$_2$, \Eqref{FN}$_1$ and \Eqref{FN2}, we receive the validity of \Eqref{BM2}$_1$. Also, from \Eqref{FN}$_2$, we infer $2r\ge q$. Finally, combining \Eqref{FN}$_1$ and \Eqref{fn}, 
we get
$$
\frac1r-\frac{2}{q}>-\frac13\,\ \ \Longrightarrow\ \ q>\frac{6r}{3+r}\,,
$$
which completes the proof of the lemma.\par\hfill$\square$\par

\renewcommand{\theequation}{3.\arabic{equation}}
\setcounter{equation}{0}
\section{Local Regularity}
The following main result guarantees local spatial regularity of distributional solutions.
\Bt Let ${\mathscr O}$ be a space-time cylinder,  and let $\bfu\in L^{2,2}(\mathscr O)$ satisfy
\be \ba{ll}\medskip
\Int{\mathscr O}{}\left\{\bfu\cdot\partial_t\bfphi+ \bfu\cdot\nabla\bfphi\cdot\bfu+ \bfu\cdot\Delta\bfphi\right\}=0\,,\ \ \mbox{for all $\bfphi\in \cald(\mathscr O)$}\,,\\
\Int{\mathscr O}{}\bfu\cdot\nabla\psi=0\,,\ \ \mbox{for all $\psi\in C_0^\infty(\mathscr O)$}\,.
\ea
\eeq{3.1}
Moreover, let $\bfu=\bfu_\sigma+\nabla\pi$ be the Helmholtz-Weyl  decomposition of $\bfu$. 
Then, if the following conditions are satisfied:
\begin{itemize}
\item[{\rm (a)}] $\bfu_\sigma\in L^{q',q}(\mathscr O)$, with
\be
\frac{2}{q'}+\frac{3}{q}=1\,,\ \ q>3\,,
\eeq{3.2}
\item[{\rm (b)}]
\be
\pi\in L^{\infty,1}(\mathscr O)\,,
\eeq{3.3}
\end{itemize}
necessarily
\be
\bfu\in L^{\infty,2}(\mathscr O')\,,\ \nabla\bfu\in L^{2,2}(\mathscr O' )\,,\ \ \mbox{for all $\mathscr O'\Subset \mathscr O$}\,. 
\eeq{3.4}
Therefore, by {\rm \cite{Ser,Stru}}, 
 $\bfu(t)$ is of class $C^\infty$ in space, for a.a. $t\in (0,T)$, and, moreover, each spatial derivative is bounded in $\bar{\mathscr O^\prime}$. 
\ET{3.1}
\Br For future reference, it is worth remarking that if $\bfu_\sigma\in L^{r',r}(\mathscr O)$ for some $r'\in [1,\infty]$, $r\in (1,\infty)$, and $\pi$ satisfies \Eqref{3.3}, by \lemmref{har}, it follows that  $\bfu\in L^{r',r}(\mathscr O')$ for arbitrary $\mathscr O'\Subset\mathscr O$. 
\ER{3.1}
\Br Let $\mathscr B:=\mathscr O\cap \real^3$. 
From \lemmref{har} and the Helmholtz-Weyl decomposition, it follows that, once we normalize $\pi$ by the request $\langle\pi(t),1\rangle_{\mathscr B}=0$, for a.a. admissible times $t$, a sufficient condition for the validity of \Eqref{3.3} is that $\bfu\in L^{\infty,s}(\mathscr O)$ for {\em some} $s>1$.
\ER{3.2}
\Br The assumption $\bfu\in L^{2,2}(\mathscr O)$ can be removed. In fact, it is enough to assume that $\bfu$ admits the Helmholtz-Weyl decomposition, with $\bfu_\sigma$ and $\pi$ satisfying (a) and \Eqref{3.3}.
\ER{3.3}
\smallskip\par
The proof of \theoref{3.1} is a straightforward consequence of the following two lemmas. 
\Bl Assume \Eqref{3.3} and that condition {\rm (a)} is satisfied for some $q\in [4,6]$. Then \Eqref{3.4} holds. 
\EL{3.1}
{\em Proof.} 
Let $Q:=B\times (t_1,t_2)\Subset \mathscr O$ be arbitrarily chosen, and let  $\bfphi$ be any element of $\cald(Q)$. For $\eta$ sufficiently small, the space-time mollifier, $\bfphi_\eta$, of $\bfphi$ belongs to $\cald(\mathscr O)$ and can be therefore replaced in \Eqref{3.1}. By a standard argument that uses the commutativity of mollification and differentiation operations, we then show that the mollifier, $\bfu_\eta$, of $\bfu$ must satisfy the following equation
\be\left.\ba{ll}\medskip 
\partial_t\bfu_{\eta}- \Delta\bfu_\eta=-\Div(\bfu\otimes\bfu)_\eta-\nabla p^{\eta}\\
\Div\bfu_\eta=0\ea\right\}
\ \ \mbox{in $Q$}\,,
\eeq{3.5}
for some smooth scalar field $p^{(\eta)}$.
Consider, next, the following linear problem
\be\ba{cc}\medskip\left.\ba{ll}\medskip
\partial_t\bfv- \Delta\bfv=-\Div(\bfu\otimes\bfu)_\eta-\nabla \phi\\
\Div\bfv=0
\ea\right\}\ \ \mbox{in $Q$\,,}\\
\bfv(x,t_1)=\0\,,\ \ x\in B\,;\ \ \bfv(x,t)=\0 \ \ \ \mbox{$(x,t)\in\partial B\times (t_1,t_2)$}\,.
\ea
\eeq{3.6}
It is well known that this problem has one (and only one) regular solution $\bfv_\eta$  (for instance, 
$\bfv_\eta\in W^{1,2}(t_1,t_2; L^2(B))\cap L^2(t_1,t_2;W^{2,2}(B)))$ which satisfies, in particular, the ``energy equation"
\be
\|\bfv_\eta\|_{L^{\infty,2}(Q)}+\|\nabla\bfv_\eta\|_{L^{2,2}(Q)}\le C\,\|(\bfu\otimes\bfu)_\eta\|_{L^{2,2}(Q)}\,.
\eeq{3.7}
Now, from \Eqref{3.5} and \Eqref{3.6}$_1$ it follows that $\bfw_\eta:=\bfu_\eta-\bfv_\eta$ solves the Stokes equations
\be\left.\ba{ll}\medskip
\partial_t\bfw_\eta-\Delta\bfw_\eta=-\nabla q^{\eta}\\
\Div\bfw_\eta=0\ea\right\}
\ \ \mbox{in $Q$}\,.
\eeq{3.8}
By \lemmref{HE} we thus deduce, in particular,
$$
\|\bfu_{\sigma\eta}\|_{L^{\infty,2}(Q')}+\|\nabla\bfu_{\sigma\eta}\|_{L^{2,2}(Q')}\le C(\|\bfu_\eta\|_{L^{1,1}(Q)}+\|\bfv_\eta\|_{L^{\infty,2}(Q)}+\|\nabla\bfv_\eta\|_{L^{2,2}(Q)})\,,  
$$
for all $Q'\Subset Q$.
We then observe that $L^{q',q}(\mathscr O)\subset L^{4,4}(\mathscr O)$, for $q\in[4,6]$. Therefore, 
since  $\|(\bfu
\otimes\bfu)_\eta\|_{L^{2,2}(Q)}\le C\,\|\bfu\|_{L^{4,4}(Q)}^2$, by assumption,  \remref{3.1}, and \Eqref{3.7}, we conclude that there exists $C_0>0$ independent of $\eta$ such that
$$
\|\bfu_{\sigma\eta}\|_{L^{\infty,2}(Q')}+\|\nabla\bfu_{\sigma\eta}\|_{L^{2,2}(Q')}\le C_0\,,
$$
from which we obtain
$$
\bfu_\sigma\in L^{\infty,2}(Q')\,,\ \nabla\bfu_\sigma \in L^{2,2}(Q')\,. 
$$
The latter, in conjunction with the assumption \Eqref{3.3} and \lemmref{har}, proves
the stated property.\par\hfill$\square$\par 
The next result is obtained by combining an argument in \cite{DiTa} with \lemmref{HE} and \lemmref{Oseen}.
\Bl Suppose \Eqref{3.3} and that condition {\rm (a)} is satisfied for some $q\in (3,4)\cup (6,\infty)$. Then  \Eqref{3.4} holds.
\EL{3.2}
{\em Proof.}  Let $Q:=B\times(t_1,t_2)$, $Q \Subset \mathscr O$ arbitrary, and let $\zeta=\zeta(x)$ be a smooth function that is 1 on arbitrarily given $B^{\prime\prime}\Subset B$ and 0 in $\real^3\backslash B$. Next, consider the vector function $\bfv$ in \Eqref{SS} with $F_{\ell j}\equiv-\zeta\,(u_\ell u_j)_\eta$. By \lemmref{Oseen}, we then deduce, on the one hand, that $\bfv\equiv\bfv_\eta$ satisfies
\be\left.\ba{ll}\medskip
\partial_t{\bfv}=\Delta\bfv+\nabla\Phi-\Div[\zeta\, (\bfu\otimes\bfu)_\eta]\\
\Div\bfv=0\ea\right\}\ \ \mbox{ in $\mathscr R_{t_1,t_2}$},
\eeq{3.9}
and also, by \Eqref{vf}, the  assumption and \remref{3.1}, that
\be
\|\bfv_\eta\|_{L^{r',r}(Q)}\le c\,\|(\bfu\otimes\bfu)_\eta\|_{L^{\frac{q'}2,\frac q2}(Q) }\le c\,\|\bfu\|_{L^{{q'},q}(Q)}^2 
\eeq{3.10}
where
\be
2r\ge q>\frac{6r}{3+r}\,,\ \ \ \ \frac{2}{r'}+\frac{3}r=1\,.
\eeq{3.11}
From \Eqref{3.5} and \Eqref{3.9} it then follows that $\bfw_\eta:=\bfu_\eta-\bfv_\eta$ satisfies \Eqref{3.8}  in $Q^{\prime\prime}=B^{\prime\prime}\times(t_1,t_2)$, where $\zeta(x)\equiv 1$. Therefore, employing \lemmref{HE}
and \Eqref{3.10}  we deduce
$$\ba{rl}\medskip
\|\bfu_{\sigma\eta}\|_{L^{r',r}(Q')}&\!\!\!\!\le c\,(\|\bfu_{\sigma\eta}\|_{L^{1,1}(Q)}+\|(\bfu\otimes\bfu)_\eta\|_{L^{\frac{q'}2,\frac q2}(Q) })\\
&\!\!\!\!\le c\,(\|\bfu\|_{L^{1,2}(Q)}+\|\bfu\|_{L^{{q'},q}(Q)}^2)\,,\ \ \mbox{arbitrary $Q'\Subset Q$.}\ea 
$$
From this relation, \Eqref{3.3} and \remref{3.1} we easily infer
\be
\|\bfu\|_{L^{r',r}(Q')}
\le c\,(\|\bfu\|_{L^{1,2}(Q)}+\|\bfu\|_{L^{{q'},q}(Q)}^2)\,,\ \ \mbox{arbitrary $Q'\Subset Q$.}
\eeq{3.12}
With equations \Eqref{3.11} and \Eqref{3.12} at our disposal, we will now show that, under the given assumptions, it is always possible to find a $\bar{q}\in [4,6]$, depending on $q$, for which condition (a) is satisfied with $\mathscr O$ replaced by a generic $\mathscr O'\Subset \mathscr O$. We will denote this modified condition (a) as (a)$^\prime$. The claimed property will then follow directly from \lemmref{3.1}. 
If $q>6$, the inequality $q>6r/(3+r)$, is satisfied for all admissible $r$. Now, suppose (a) holds for $q\equiv q_0>6$. Then,   from \Eqref{3.11} and \Eqref{3.12}  it follows that (a)$^\prime$ holds also for any $q\ge q_1:=q_0/2$. If $q_1\le 6$, we are done, otherwise we use \Eqref{3.12} with $q=q_1$ and derive, again from \Eqref{3.11}, that (a)$^\prime$ holds for any $q\ge q_2:=q_1/2=q_0/2^2$. If $q_2\le 6$ the proof is complete; if not, we iterate this process a finite number of times, until we find an integer $\bar{n}$ such that condition (a)$^\prime$ is satisfied for all $q\ge q_{\bar{n}}:=q_0/2^{\bar{n}}\le 6$ and, as a result, by some $\bar{q}\in[4,6]$. Note that, in this process, the domains $Q$ and $Q'$ may ``shrink" at each step by an arbitrarily given small amount, and the constant $c$ in \Eqref{3.12} may also change. However, $Q'$ and $Q$ are themselves arbitrary and the number of steps is finite. Therefore,  the proof in the case $q>6$ is completed. Next, we take $q\in(3,4)$ and, to show the result, we employ the elegant argument used in \cite[p. 18-19]{DiTa}. Precisely, set $q:=q_k, q'_k:=q'$, and $q_{k-1}:=r, q^\prime_{k-1}=r^\prime$, with $r,r^\prime$ satisfying \Eqref{3.11}--\Eqref{3.12}, namely,
\be      
4>q_k>\frac{6q_{k-1}}{3+q_{k-1}}\,,\ \ \frac{2}{q'_{k-1}}+\frac 3{q_{k-1}}=1\,,
\eeq{3.13}
and
\be
\|\bfu\|_{L^{q'_{k-1},q_{k-1}}(Q')}
\le c\,(\|\bfu\|_{L^{1,2}(Q)}+\|\bfu\|_{L^{{q'_k},q_k}(Q)}^2)\,.
\eeq{3.14}
Consider the sequence 
\be
c_k=\frac{6c_{k-1}}{3+c_{k-1}}\,,\ \ c_0=4\,,
\eeq{3.15}
and observe that $\lim_{k\to\infty}c_k=3$.
Thus, from  \Eqref{3.13} and \Eqref{3.15} it follows that, given $i\in \mathcal I:=\{0,\ldots,k-1\}$, for any $q_{k-i}>c_{k-i}$ there is a $q_{k-i-1}>c_{k-i-1}$ such that the pair $(q_{k-i-1}^\prime,q_{k-i-1})$ satisfies \Eqref{3.13}$_2$ and also \Eqref{3.14}. Therefore, varying $i$ in $\cali$, after $k$ iteration we conclude that there is $q_0>4$ such that
$$
\|\bfu\|_{L^{q'_{0},q_{0}}(Q')}
\le c\,(\|\bfu\|_{L^{1,2}(Q)}+\|\bfu\|_{L^{{q'_k},q_k}(Q)}^2)\,.
$$
with $q'_0, q_0$ satisfying \Eqref{3.13}$_2$. In this iterative scheme,  $Q'$ may ``shrink" at each step by an arbitrarily given small amount as well as the constant $c$ in \Eqref{3.14} may change but,  as previously mentioned, this does not affect the conclusion of the lemma which is thus completely proven.\par\hfill$\square$\par
We conclude this note with some remarks concerning  {\em time regularity} of distributional solutions.
Although, as mentioned in the Introduction, the {time regularity} of $\bfu$ cannot be established without making appropriate assumptions about it for $\nabla\pi:={\sf Q}\,\bfu$, we can nevertheless ask whether something can be said about the time regularity of $\bfu_\sigma$.The answer is yes. In fact,  under the assumptions of \theoref{3.1}, we obtain
$$
\bfu\cdot\nabla\bfu\in L^{s,s'}(Q)\,,\ \ \mbox{for all $s,s'\in [1,\infty]$\,,\ \ all $Q\Subset \mathscr O$.}
$$
Therefore, from \cite[Theorem 2.8]{GiSo}, we deduce that problem \Eqref{3.6} has one (and only one) solution
$$
\bfv\in W^{1,s^\prime}(L^s)\cap L^{s^\prime}(W^{2,s})\,,\ \ s^\prime,s\in (1,\infty)\,,
$$
that also depends continuously on $\|\bfu\cdot\nabla\bfu\|_{L^{s^\prime,s}}$. We may then use the same argument employed in the proof of \lemmref{3.1} and show that
\be
\partial_t\bfu_\sigma\in L^{s^\prime,s}(Q')\ \ \mbox{for all $Q^\prime\Subset Q$}\,. 
\eeq{3.16}
Notice that if, in addition, we suppose $\partial_t\pi\in L^{s',1}$,  then, from \lemmref{har} (with $\partial_t\pi$ in place of $\pi$) and \Eqref{3.16}, we conclude  $\partial_t\bfu\in L^{s^\prime,s}$. This fact suggests the following {\em conjecture}: Under the assumptions of \theoref{3.1}, suppose, further,  $\partial_t^k\pi\in L^{r,1}$, $k=1,\ldots m$, $m>1$. Then $\partial_t^k\bfu\in L^{r,s}$, for all $s>1$.   
\medskip\par\noindent
{\bf Acknowledgment.} Work  partially supported by National Science Foundation (US) Grant DMS--2307811.

\ed